\documentclass[reqno,12pt]{amsart}
\usepackage{indentfirst}
\usepackage{lmodern}
\usepackage{amsmath}
\usepackage[T1]{fontenc}
\usepackage{amssymb}			
\usepackage{amsthm}
\usepackage{amsfonts}
\usepackage{enumerate}
\usepackage{enumitem}
\usepackage{setspace}
\usepackage{mathtools}
\usepackage[body={7in,9.5in},top=1in, left=0.8in ]{geometry}
\newcommand{\bn}{\mathbb{N}}

\newcommand{\ov}[1]{\overline{#1}}

\newcommand{\pref}[1]{\textup{(\ref{#1})}}
\newcommand{\ord}[1]{\textup{ord}({#1})}

\newtheorem{thrm}{Theorem}

\newtheorem{thm}[thrm]{Theorem}
\newtheorem{definition}[thrm]{Definition}
\newtheorem{remark}[thrm]{Remark}

\newtheorem{df}[thrm]{Definition}		

\newtheorem{rmk}[thrm]{Remark}

\newtheorem{prop}[thrm]{Proposition}
\newtheorem{exam}[thrm]{Example}

\usepackage[nobysame,initials]{amsrefs}

\newcommand{\X}{{\mathbb{X}}}

\newcommand{\K}{{\mathbb{K}}}
\newcommand{\R}{{\mathbb{R}}}
\newcommand{\C}{{\mathbb{C}}}
\newcommand{\Z}{{\mathbb{Z}}}
\newcommand{\N}{{\mathbb{N}}}

\title[Multiplicative roots of formal power series]{A note on multiplicative roots of multivariable formal power series}

\author[P. Ma\'ckowiak]{Piotr Ma\'ckowiak}
\author[M. Mokatren]{Motaz Mokatren}

\address[P. Ma\'ckowiak]{Department of Nonlinear Analysis and Applied Topology\\
  Faculty of Mathematics and Computer Science\\
  Adam Mickiewicz University, Pozna\'n\\
  Uniwersytetu Pozna\'nskiego~4\\
  61-614 Pozna\'n\\
  Poland}

\address[M. Mokatren]{Department of Electrical and Electronics Engineering\\
	Kinneret College\\ 
	Zemach, Emek HaYarden Mobile Post 15132\\
	Israel}

\subjclass[2010]{Primary: 13F25; Secondary: 13J05}

\keywords{multivariable formal power series, formal root series, multiplicative root}

\email[P.~Ma\'ckowiak]{piotr.mackowiak@amu.edu.pl}
\email[M. Mokatren]{motaz\_mok@mx.kinneret.ac.il}

\begin{document}

\maketitle
\begin{abstract}
Suppose that we are given a formal power series of many variables with coefficients in $\R$ (or $\C$) and we want to compute its $n$-th (multiplicative) root. As can be expected coefficients of the root have to satisfy a system of infinitely many equations. We present such a system of equations that strictly corresponds with the system for $n$-th of a formal power series of one variable. With help of an example we show that the case of formal power series of many variables is very different from the one variable case with respect to the existence of roots.
\end{abstract}

\section{Introduction}
In this paper we are going to address the following simple question: for a given formal power series of many variables $f$, is it possible to determine its multiplicative $n$-root, that is, to find a formal power series $g$ for which $g^n=f$ or to determine that the series has no such a root? 
The case of one-variable is fully described in \cite{XD}. It turns out that in this case the existence of $n$-th root fully depends on the existence of $n$-th root of the first nonzero term. Our main goal is to check whether similar characterization holds for multivariable formal power series.

W.G.Brown has formulated a theorem on the existence of a square root of a power series of two variables \cite{WB65}, which he later extended to the existence of an $n$-th root \cite{WB67}. Brown's theorems have thoroughly examined the existence of roots for both unit series (nonzero constant term) and nonunit series (zero constant term). According to theses theorems a formal power series $Q$ has an $n$-th root if and only if it admits a factorization of the form $Q=P^nR$, where $R$ is a unit formal power series possessing $n$-th root (see Thm 2.2 of \cite{WB67}). However, in general it is not clear how to determine $P^n$ and/or $R$ for a given $Q$.

In \cite{N69} I.Niven popularized the theory of formal power series in one variable. Among the others, he briefly discussed how to derive the coefficients of an $n$-th root of some specific unit formal power series. As we already mentioned, Gan and Bugajewski provided an algorithm for the computation of coefficients of $n$-th root of formal power series (unit or not), assuming such a root exists. In the paper \cite{XD} there are also given criteria for the existence of roots of one-variable formal power series.  

It may seem strange but as we are aware of there are no multivariable case generalizations of the results from \cite{XD} in the literature with the exception that in \cite{Hunk2019} one can find a result on the existence and computation of $n$-th root of a unit formal power series whose constant term is $1$.

Let us also notice that multiplicative roots of formal power series are of some interest. Applications of multiplicative roots can be found in enumerative graph theory (see e.g \cite{BTW} or \cite{TW}), existence theorems for solutions of quadratic equations (see \cite{WB65} or \cite{TW}) or Riordan group theory (see \cite{CM}). 

The paper is organized as follows. The next section sets notation and gives basic definitions. Section 3 presents formulas for the computation of a root of a unit multivariable formal power series and contains a necessary and sufficient condition for the existence of an $n$-th root of multivariable formal power series. At the end of that section an example is presented. The example shows that one variable formal power series and many variable formal power series behave differently with respect to the existence on an $n$-th root.
\section{Preliminaries}
The goal of this section is to recall and introduce basic definitions and facts which will be needed in the sequel. The notation is more-or-less borrowed from \cite{BBXM}.

Symbol $\Z$ denotes the set of integer numbers, $\N$ is the set of positive integers and $\N_0:=\N\cup\{0\}$. If $n\in \N$, then $[n]:=\{1,\ldots,n\}$ and $[n]_0:=[n]\cup \{0\}$. 
For any set $S$ by $(S)^n$ or $S^n$ we denote the cartesian product of $n$ copies of $S$, $n\in \bn$; $S^0:=\emptyset$.

Let us fix a $q\in \N$ and denote by $C_m$ the set of all nonnegative
integer solutions $c_1,\ldots,c_q$ to the equation $c_1+\ldots+c_q=m$ for $m\in \N_0$, that is
\begin{equation}\label{eqn:exponents}
C_m:=\{c=(c_1,\ldots, c_q)\in \N_0^q:\, c_1+\ldots+c_q=m\},
\end{equation}
where $\N_0^q:=(\N_0)^q$. The number $r(m):=\binom{m+q-1}{q-1}$ is the number of elements of $C_m$ (see \cite[p. 25]{S}). Hence, we can write
$$C_m=\{c^m_1,\ldots, c^m_{r(m)}\}\subset \bn^q_0.$$
It is clear that $r(0)=1$ and $C_m\cap C_{m'}=\emptyset $ for $m\neq m'$.
Obviously, $C:=\N^q_0=\bigcup_{m\in \bn_0}C_m$ and for each $c\in C$ there is exactly one $m\in \bn_0$ for which $c\in C_{m}$ and we denote it by $m(c)$. It is clear that $m(c)=c_1+\ldots+c_q$. For any $m\in \bn_0$, we set $C_{[m]}:=\bigcup_{i\in [m]_0}C_i$. We also define $\theta:=(0,...,0)\in \N^q_0$.

If $A$ is an $m\times n$ matrix, the element from $i$-th row and $j$-th column of $A$ is denoted as $a_{ij}$; $A^T$ is the transpose of $A$.

We take on the convention that $\sum_{x\in \emptyset}:=0$ and $\inf \emptyset:=+\infty$. Moreover, $\infty+\infty:=\infty$ and $n+\infty)=\infty+n=\infty$ for all $n\in \N_0$. 
We say that $a=(a_1,\ldots,a_q)\in C$ is lexicographically less (greater) than $b=(b_1,\ldots,b_q)\in C$, $a\neq b$, if $a_i<b_i$ ($a_i>b_i$) for the first coordinate $i\in [q]$ for which $a_i\neq  b_i$.
 
Let $\K$ stand for the field $\R$ of real (or $\C$ of complex) numbers.

\begin{definition}
A~formal power series $f$ of $q$-variables ($q$-fps) over $\K$ is a function acting from $\bn^q_0$ to $\K$. The value $f_c:=f(c)\in \K$ is called the $c$-th coefficient of $f$, $c\in C$. The set of all $q$-fps (over $\K$) is denoted as $\X_q$. 
\end{definition}
It is clear that $q$-fps defined above, if we do not consider any topological issues connected with formal power series, is a standard multivariable formal power series over $\K$ (cf. \cite[Definition 2.1]{Hunk2019} or \cite[Definition 8.1]{Samb}). Let $X^c:=x_1^{c_1}\ldots x_q^{c_q},\, x:=(x_1,\ldots,x_q)$, where $x_i$ stands for a formal variable 
$i\in [q]$, $c\in C$. Using summation notation, $q$-fps $f$ can be written as the formal sum $$f(x):=\sum_{c\in C}f_cX^c.$$
Observe that monomials of the form $aX^c$, where $a\in \K$ and $c\in C$, are $q$-fps whose all coefficients, except perhaps the $c$-th one, are $0$. In particular, value of the only nonzero coefficient of the monomial $X^c$ is $1$.

\begin{rmk}
Since $q$-fps are functions whose codomain is $\K$ we can add or subtract $q$-fps pointwise. Similarly we can multiply any $q$-fps by a scalar $\lambda$ to get $\lambda f$ according to the formula $(\lambda f)_c:=\lambda f_c, c\in C$. For basic properties of multivariable formal power series see e.g. \cite{Hunk2019}. 
\end{rmk}
\begin{df}
Let $f\in \X_q$. $f$ is called a unit if $f_\theta\neq 0$, otherwise $f$ is called a nonunit. We define the order of $f$ by $\ord{f}:=\inf\{k\in \N_0:\, \textrm{ there exists }c\in C_k\textrm{ for which }f_c\neq 0\}.$ 
\end{df}
Observe that if $\ord{f}=0$, the initial block of $f$ is of the form $f_\theta X^\theta$.

The following definitions are vital for our paper.
\begin{df}
Let $f\in \X_q$, $\ord{f}\in \N_0$. The initial block of $f$ is $q$-fps defined as $\sum_{c\in C_{\ord{f}}} f_c X^c$.
\end{df}
\begin{df}\label{df:prod}
For $q$-fps $f(x)=\sum_{c\in C}f_cX^c, \, g(x)=\sum_{c\in C}g_cX^c\in \X_q$, $q\in \bn$, the $q$-dimensional Cauchy product (in short: product) of $f$ and $g$ is a $q$-fps $h=fg$ defined as
$$h(x):=\sum_{c\in C}\underbrace{\left(\sum_{a,\,b\in C:\,a+b=c}f_ag_b\right)}_{h_c:=}X^c.$$
\end{df}

\begin{df}\label{df:root}
Let $f\in \X_q$ and $n\in \bn$. Any $g\in \X_q$, for which $g^n=f$, where $g^n:=\underbrace{g\ldots g}_{n\times}$, is called an $n$-th root of $f$. 
\end{df}
It is clear that $g$ is $n$-th root of $f$ provided that $n$-th power of $g$ equals $f$: for any $c\in C$ it holds $f_c=(g^n)_c$.

We shall need the following simple counterpart of a well--known property of $1$-fps. 
\begin{prop}\label{prop:ord}
Let $f, g\in \X_q$. Then $\ord{fg}=\ord{f}+\ord{g}$.
\end{prop}
\begin{proof}
If $\ord f=\infty$ or $\ord g=\infty$, then the claim is true. The same is true if at least one order is $0$. Suppose that $\ord{f}=m\in\bn$, $\ord{g}=n\in \N$. Let $c\in C_k$, $k<m+n$, $k\in \N_0$. By definition of the Cauchy product $(fg)_c=\sum_{a+b=c}f_ag_b$, where $a, b\in C$ (recall that $C=\N^q_0$). Obviously $a\in \C_{m(a)}$ and $b\in C_{m(b)}$ and $m(a)+m(b)=m(c)<k$. Hence, $f_a=0,\,g_b=0$ for such $a$ and $b$. Therefore $\ord{fg}\geq m+n$, provided $fg\neq \theta$, where $\theta$ represents a $q$-fps whose all coefficients equal $0$. Since the orders of $f$ and $g$ are $m$ and $n$, respectively, then there exist $c'\in C_{m}$ and $c''\in C_{n}$ such that $f_{c'}\neq 0\neq g_{c''}$. Without loss of generality we may assume that $c'$ and $c''$ are lexicographically the least ones in $C_m$ and $C_n$, respectively, for which $f_{c'}\neq 0$ and $g_{c''}\neq 0$. Let $c:=c'+c''$. Then $c\in C_{m+n}$. We also have $(fg)_c=\sum_{a+b=c}f_ag_b=\sum_{a+b=c,\, a\in C_m, b\in C_m}f_{a}g_{b}$, since if $m(a)<m$, then $f_a=0$ (and similarly for $m(b)<n$, $g_{b}=0$). But if $a\in C_m$, $b\in C_n$, $a+b=c$ and $a$ ($b$) is lexicographically greater than $c'$, then $b$ ($a$) is lexicographically less than $c''$ which implies that in such a case $f_ag_b=0$. Thus, $(fg)_c=f_{c'}g_{c''}$ and the claim follows.
\end{proof}
\begin{remark}It is a straightforward consequence of Proposition \ref{prop:ord} that the ring of formal power series (with respect to addition and multiplication) is an integral domain (cf. \cite{Hunk2019}).
\end{remark}
\section{Roots}
\subsection{Units}
If $f$ is $1$-fps there is no need to discern between a vector of formal variables $x$ and the formal product of powers of formal variables $X^c$. We identify any $1$-fps $f$ with $f(x)=\sum_{i\in\N_0}f_ix^i$. By Theorems 1 and 2 from \cite{XD} we can state the following theorem characterizing the existence of roots of one variable formal power series.
\begin{thm}\label{thm:BGfps}
Let $f\in \X_1$, $\ord{f}\in\N_0$, and $n\in \N$. There exists an $n$-th root of $f$ if and only if the its initial block, $h(x):=f_{\ord{f}}x^{\ord{f}}$, possesses an $n$-th root.  
\end{thm}
According to Theorem \ref{thm:BGfps}, whether an $n$-th root of an $1$-fps $f$ exists or not can be decided by examining the first (in  natural ordering of $\N_0$) nonzero term of $f$, that is, $f_{\ord{f}}$. Observe, that $h(x)$ is the minimal monomial of $f$. Hence,  $f^{1/n}$ exists if and only if the minimal monomial of $f$ possesses an $n$-th root.

It is not clear what ordering should be treated as a natural ordering of elements of $C=\N^q_0$ if $q>1$. However, let us observe that if we want to compute $c$-th and $c'$-th terms of the product $fg$, where $m(c)<m(c')$, then to calculate $(fg)_{c'}$ we need more information and operations than for $(fg)_c$. That's why we think that a kind of natural partial ordering of elements of $C$ is the following:  for $c,c'\in C$, $c$ precedes $c'$ if the sum of components of $c$ is less than or equal to the sum of components of $c'$, that is, $m(c)\leq m(c')$. If we additionally impose the lexicographic ordering of elements of $C_m$, $m\in \N_0$, then we get the graded lexicographic ordering of $C$ (see \cite{CLO}*{Ch.2 $\mathsection$2}) which is a linear order on $C$.

If a $q$-fps has nonzero constant term, the situation is very similar to the case of one variable (cf. \cite{XD}) with the exception that now we examine initial blocks rather than initial terms.
\begin{thm}\label{thm:qunits}
Let  $f(x)=\sum_{c\in C}f_c X^{c}$ be a unit $q$-fps, $q\in \N$. Let also $n\in \N$ be given. An $n$-th root of $f$ exists if and only if $h(x):=f_\theta X^\theta$, the initial block of $f$, possesses an $n$-th root. 
\end{thm}
\begin{proof}
Let $f(x)=\sum_{c\in C}f_c X^{c}$, $f_{\theta}\neq{0}$, possess an $n$-th root $g=f^{1/n}$. Thus, by definition of the Cauchy product, $f_\theta=\underbrace{g_\theta\ldots g_\theta}_{n\times}$, which implies that $h(x)=(g_\theta X^\theta)^n$.

Now, let us assume that there is $t\in \K$ such that $t^n=f_\theta$, so $h(x)=f_\theta X^{\theta}$ has an $n$-th root $tX^\theta$. We shall now construct a $q$-fps $g\in \X_q$ for which $g^n=f$. To this goal let 
\begin{equation}\label{eq:u1}
g_\theta:=t.
\end{equation} 
Hence $(g^n)_\theta=t^n=f_\theta$. Let $c\in C_1$. Observe that $(g^n)_c=f_c$ if and only if $\sum g_{c^1}\ldots g_{c^n}=f_c$, where the summation extends over all solutions $(c^1,\ldots,c^n)\in C^n$ of the equation ${c^1+\ldots+c^n=c}$. But we have for the sum that $\sum g_{c^1}\ldots g_{c^n}=ng^{n-1}_\theta g_{c^1}$, so it suffices to define 
\begin{equation}\label{eq:u2}
g_{c}:=f_c(ng^{n-1}_\theta)^{-1}.
\end{equation} 
Suppose that we have defined $g_c$ for all $c\in C_{[m]}$, so that we have $f_c=(g^n)_c$, $c\in C_{[m]}$. We shall define $g_c$, $c\in C_{m+1}$, in such a way that $f_c=(g^n)_c$ for $c\in C_{[m+1]}$. This will complete the proof. Let us fix $c\in C_{m+1}$. We have that $f_c=(g^n)_c$ if and only if $\sum g_{c^1}\ldots g_{c^n}=f_c$, where the summation extends over all solutions $(c^1,\ldots,c^n)\in C^n$ of the equation ${c^1+\ldots+c^n=c}$. But $\sum g_{c^1}\ldots g_{c^n}=ng^{n-1}_\theta g_c+S_m$, where $S_m:=\sum_{(c^1,\ldots,c^n)\in (C_{[m]})^n:\, c^1+\ldots+c^n=c} g_{c^1}\ldots g_{c^n}$. Therefore it suffices to set 
\begin{equation}\label{eq:u3}g_c:=(f_c-S_m)(ng^{n-1}_\theta)^{-1}.\end{equation}
\end{proof}
Observe that equations \pref{eq:u1}-\pref{eq:u3} allow for the computation of all coefficients of $f^{1/n}$ under the given initial value $t$.
\subsection{Nonunits}
Let us start with the following theorem which we state without proof as it is a simple consequence of the definition of a root, associativity of multiplication and Theorem \ref{thm:qunits}.
\begin{thm}\label{thm:simple}
Let $f(x)=\sum_{c\in C}f_c X^{c}\in \X_q$ with $\ord{f}=m\in \N$ and  let $n\in \bn$. Suppose that $f$ has initial block of the form $f_{\alpha}X^\alpha$, for some $\alpha=(\alpha_1,\ldots,\alpha_q)\in C_m$ (see formula \pref{eqn:exponents}), and there exists $g(x)=\sum_{c\in C}g_c X^{c}\in \X_q$ such that $g_\theta=1$ and $f(x)=(f_\alpha X^\alpha)g(x)$. 

Then $f$ possesses an $n$-th root if and only if 
$f_\alpha X^\alpha$ possesses an $n$-th root.
Moreover, if $h(x)=\sum_{c\in C}h_cX^c$ is an $n$-th root of $f$, then $n\times \ord{h}=m$, $h$ has the initial block of the form $h_{\ov{c}} X^{\ov{c}}$, for some $\ov{c}\in C_{\ord{h}}$, and we have $(h_{\ov{c}})^n=f_\alpha$ and $n\ov{c}=\alpha$. And conversely, if $h\in \X_q$ , $n\times \ord{h}=m$, $h$ has the initial block of the form $h_{\ov{c}} X^{\ov{c}}$, for some $\ov{c}\in C_{\ord{h}}$, and we have $(h_{\ov{c}})^n=f_\alpha$ and $n\ov{c}=\alpha$, then $f=h^n$. 
\end{thm}

As we already mentioned in the introduction, Brown's theorem examines the existence of roots for nonunit power series. It is asserted there that a formal power series $Q$ has an $n$-th root if and only if it can be factored as $Q= P^k R$ \cite[Theorem 2.2]{WB67}, where $R_\theta=1$ and $P$ is a formal power series. However, in general it is not clear how to determine the formal power series $P$ (or $R$) appearing in the decomposition. To see this, just consider the expanded form of $n$-th power, for some large $n$, of a formal power series $\sum f_cX^c$, where no coefficient $f_c$ vanishes and there is no clear regularity of coefficients $f_c$, $c\in C$. Thus we are given $f^n$ in an expanded form, and we know that its $n$-th root exists (by construction). But if we did not know how the fps was produced, we would not be able to easily find its $n$-th root or determine if an $n$-th root exists. Observe also that our simple Theorem \ref{thm:simple} is not applicable here. In what follows we will try to tackle the problem of finding a root and/or determining its existence in a constructive way. As it will turn out, in general, the situation is rather hopeless.

Let $f,\, g\in \X_q$ and $n\in \N$ be fixed. By Proposition \ref{prop:ord} it follows that if $f^{1/n}=g$, that is, $g^n=f$, then $n\times \ord{g}=\ord{f}$. Thus, divisibility of $\ord{f}$ by $n$ and $\ord{g}$ is a necessary condition for $g$ to be an $n$-th root of $f$. 

Denote $s:=\ord{f},\, m:=\ord{g}$, and assume that $mn=s>0$. Suppose that $g^n=f$, that is, $g$ is an $n$-th order multiplicative root of $f$. By definition of the (Cauchy) multiplication of $q$-fps, for $c\in C$,
\begin{equation}\label{eq:1}
f_c=\sum\{g_{c^1}\ldots g_{c^n}:\, (c^1,\ldots,c^n)\in C^n \textrm{ and }c^1+\ldots+c^n=c\}.
\end{equation}
In particular, for $c\in C_{mn}$, for any solution $(c^1,\ldots,c^n)\in C^n$ of the equation $c^1+\ldots+c^n=c$ such that, for some $i\in [n]$, $c^i\in C_k$, where $k<m$, we have $g_{c^1}\ldots g_{c^n}=0$. This implies
that, for $c\in C_{mn}$, the only components in the sum on the right-hand side of \pref{eq:1} that may contribute a nonzero value to $f_c$ are those that correspond to solutions of $c^1+\ldots+c^n=c$, where $c^i\in
C_m\cup C_{m+1}\cup\ldots,\, i\in [n]$. On the other hand, since $c\in C_{mn}$ and $c=\underbrace{c^1+\ldots+c^n}_{n \times}$, for the just mentioned possibly nonzero components it holds $c^i\in C_m$ for every $i\in [n]$.
Thus we have, for $c\in C_{mn}$,
\begin{equation}\label{eq:1'}
f_c=\sum_{\scriptsize\begin{array}{c} (c^1,\ldots,c^n)\in (C_{m})^n:\\\sum_{i\in [n]}c^i=c\end{array}} g_{c^1}\ldots g_{c^n}.
\end{equation}
Let us now consider $c\in C_{mn+1}$. Hence, $f_c$ is a sum of the form \pref{eq:1}. Since $c\in C_{mn+1}$ and there are $n$ (not necessarily different) elements $c^1,\ldots, c^n\in C$ such that
$c^1+\ldots+c^n=c$, the product $g_{c^1}\ldots g_{c^n}$ may be nonzero only if exactly $n-1$ out of the $n$ elements $c^1,\ldots,c^n$ belong to $C_m$ and the remaining element is a member of $C_{m+1}$. This is a consequence
of the fact that if at least one $c^i$ belongs to $C_k$ with $k>m+1$, then some component $c^j$ in the sum $c^1+\ldots+c^n(=c)$ belongs to $C_l$ with $l<m$ which entails that $g_{c^j}=0$. Therefore we can write
\begin{multline*}f_c=
\sum_{j=1}^{r(m+1)}\sum_{\scriptsize\begin{array}{c}(c^1,\ldots,c^n)\in (C_{[m+1]})^n:\\\sum_{i\in [n]}c^i=c,\\\exists i\in [n]: \, c^i=c^{m+1}_j\end{array}}g_{c^1}\ldots g_{c^n}=\sum_{j=1}^{r(m+1)}g_{c^{m+1}_j}\sum_{\scriptsize\begin{array}{c} (c^1,\ldots,c^n)\in (C_{m})^{n-1}:\\\sum_{i\in [n-1]}c^i=c-c^{m+1}_j\end{array}}ng_{c^1}\ldots g_{c^{n-1}}.\end{multline*}
Hence, for $c\in C_{mn+1}$,
$$f_c=H_c^1g_{c^{m+1}_1}+\ldots+H_c^{r(m+1)}g_{c^{m+1}_{r(m+1)}},$$ where $H_c^j:=n\sum g_{c^1}\ldots g_{c^{n-1}},$ and the sum extends over solutions $c^1,\ldots,c^{n-1}\in C_m$ of the equation
$c^1+\ldots+c^{n-1}=c-c^{m+1}_j,\, j\in [r(m+1)]$. Since $C_{mn+1}=\{c^{mn+1}_1,\ldots,c^{mn+1}_{r(mn+1)}\}$, using matrix notation we obtain
\begin{equation}\label{eq:2}[f]_{mn+1}=H^{m}[g]_{m+1},\end{equation}
where $[f]_{mn+1}=[f_{c^{mn+1}_1}\ldots f_{c^{mn+1}_{r(mn+1)}}]^T$, $[g]_{m+1}=[g_{c^{m+1}_1}\ldots g_{c^{m+1}_{r(m+1)}}]^T$, $H^{m}$ denotes the $r(mn+1)\times r(m+1)$ matrix with elements $h^{m}_{ij}:=H_{c^{mn+1}_i}^j$,
$i\in [r(mn+1)]$, $j\in [r(m+1)]$.



Let us now consider $c\in C_{mn+k}$, $k\geq 2$. We have
\begin{multline*}f_c=\sum_{\scriptsize\begin{array}{c}(c^1,\ldots,c^n)\in (C_{[m+k]})^n:\\\sum_{i\in [n]}c^i=c\end{array}}g_{c^1}\ldots g_{c^n}\\=\sum_{\scriptsize\begin{array}{c}(c^1,\ldots,c^n)\in (C_{m}\cup\ldots \cup C_{m+k-1})^n:\\\sum_{i\in
[n]}c^i=c\end{array}}g_{c^1}\ldots g_{c^n}+\sum_{\scriptsize\begin{array}{c}(c^1,\ldots,c^n)\in (C_{m}\cup C_{k})^n:\\\sum_{i\in [n]}c^i=c\end{array}}g_{c^1}\ldots g_{c^n}\\=\underbrace{\sum_{\scriptsize\begin{array}{c}(c^1,\ldots,c^n)\in (C_{m}\cup\ldots \cup C_{m+k-1})^n:\\\sum_{i\in [n]}c^i=c\end{array}}g_{c^1}\ldots g_{c^n}}_{G^{m+k-1}_c:=}+\sum_{j=1}^{r(m+k)}g_{c^{m+k}_j}\underbrace{\sum_{\scriptsize\begin{array}{c}(c^1,\ldots,c^{n-1})\in (C_{m})^{n-1}:\\\sum_{i\in [n-1]}c^i=c-c^{m+k}_j\end{array}}ng_{c^1}\ldots
g_{c^{n-1}}}_{H^j_c:=}.\end{multline*}

Setting $[f]_{mn+k}:=[f_{c^{mn+k}_1}\ldots f_{c^{mn+k}_{r(mn+k+1)}}]^T$, $[g]_{m+k}:=[g_{c^{m+k}_1}\ldots g_{c^{m+k}_{r(m+k)}}]^T$, 
$G^{m+k-1}:=[G^{m+k-1}_{c^{mn+k}_1}\ldots G^{m+k-1}_{c^{mn+k}_{r(mn+k)}}]^T$ and defining
entries of an $r(mn+k)\times r(m+k)$ matrix $H^{m+k-1}$ by $h^{m+k-1}_{ij}:=H_{c^{mn+k}_i}^j$, $i\in [r(mn+k)]$, $j\in [r(m+k)]$, we get
\begin{equation}\label{eq:4}[f]_{mn+k}=G^{m+k-1}+H^{m+k-1}[g]_{m+k}.\end{equation}

The right--hand sides of equations \pref{eq:1'}, \pref{eq:2} and \pref{eq:4}, $k\in\N$, are $c$-th coefficients of $n$-th power of $g$, so we conclude the above discussion with the following theorem.
\begin{thm}\label{thm:ns}
Let $f,\,g\in \X_q$, $s:=\ord{f}>0$, $m:=\ord{g}$, $n\in \bn$, and $s=mn$. Then $g$ is an $n$-th root of $f$ if and only if equations \pref{eq:1'}, \pref{eq:2} and \pref{eq:4}, $k\in\N$, are satisfied.
\end{thm}
Notice that coefficients $f_c$, $c\in C_{mn}$, appearing in equation \pref{eq:1'} are coefficients of the initial block of $f$. Similarly, $g_a, \, a\in C_m$, that appear on the right-hand side in equation \pref{eq:1'} are coefficients defining the initial block of $g$.

Let us now turn to the question how to calculate an $n$-th root of a nonunit $q$-fps $f$ or determine that such a root does not exist.
Observe that matrices $H^{m}$ and $H^{m+k-1}$ in equations \pref{eq:2} and \pref{eq:4}, respectively, depend only on coefficients $f_c,\, c\in C_m$, while the column vector $G^{m+k-1}$ in equation \pref{eq:4} depends on terms $f_c,\, c\in C_m\cup\ldots\cup C_{m+k-1}$. Since all matrices in equations \pref{eq:2} and \pref{eq:4} depend on coefficients of the initial block of $g$, we need to have these coefficients first so that we can compute coefficients $g_c$ for $c\in C_k,\,k\geq m+1$. The coefficients $g_c$, $c\in C_m$, can be obtained by solving equations \pref{eq:1'}, $c\in C_{mn}$. Notice that right--hand side of equation \pref{eq:1'} is a polynomial equation in variables $g_a,\, a\in C_{m}$, and the solution set of this equation forms an affine variety (we refer the reader to the monograph \cite{CLO} for more information on affine varietes). Suppose that we are given a solution to equation \pref{eq:1'}. Then we can use \pref{eq:2} to determine $g_a,\, a\in C_{m+1}$ (so that \pref{eq:2} be satisfied for all $c\in C_m$). Notice that there is a linear dependence between $[g]_{m+1}$ and $[f]_{mn+1}$ and it is possible that there are infinitely many solutions $[g]_{m+1}$ to \pref{eq:2}. Then we can try to compute the coefficients $[g]_k$, $k\geq m+1$, using equation \pref{eq:4}. Again, it is possible that at some stage of computations there are infinitely many solutions of \pref{eq:4} for $[g]_k$.

Formulas \pref{eq:1'}, \pref{eq:2} and \pref{eq:4} allow us to prove the following theorem.
\begin{thm}\label{thm:necc}
Let $f\in \X_q$, $s:=\ord{f}>0$, $n\in \bn$, and $n$ divides $s$. If $f$ has an $n$-th root, the initial block of $f$, $h(x):=\sum_{c\in C_{\ord{f}}} f_c X^c$, possesses an $n$-th root, as well.
\end{thm}
\begin{proof}
Let $g\in \X_q$ be an $n$-th root of $f$. Then $r(x):=\sum_{c\in C_{\ord{g}}} g_c X^c$ is an $n$-th root of $h$. To see this observe that $\ord{h}=n\times \ord{r}$ and that equation \pref{eq:1'} is satisfied (with $f$ replaced with $h$ and $g$ replaced with $r$). Moreover, equations \pref{eq:2} and \pref{eq:4} are also satisfied (with the mentioned replacements), since we have $0$'s on both sides of every equation. 
\end{proof}

We conclude the paper with the following example that shows that it is impossible to reverse the implication of Theorem \ref{thm:necc}. 
This example reveals the main difference between $1$-fps and $q$-fps, $q\geq 2$, with respect to the existence of an $n$-th root of $f$ and the existence of $n$-th root of the initial block of $f$.

\begin{exam}\label{xmpl:1}
Let $g(x,y):=\underbrace{(xy-\frac{x^3y^3}{3!}+\frac{x^5y^5}{5!}-\frac{x^7y^7}{7!}+\ldots)}_{h(x,y):=}+x^k\in \X_2$, where $k\geq 2$ is fixed. Here we identify $x^k$ with $x^ky^0$ and $x^1$ $(y^1)$ with $x$ $(y)$. Define $f(x,y):=g^2(x,y)$. Then  
$$f(x,y)=h^2(x,y)+2x^kh(x,y)+x^{2k},$$
and $\ord{f}=4$. By the construction it follows that $g$, whose order is $\ord{g}=2$, is a square root of $f$. Let $$f_q(x,y):=h^2(x,y)+2x^kh(x,y)+x^q,$$ where $q>2k$ is an odd number. 
Notice that $\ord{h^2}=4$ and $\ord{2x^kh(x,y)}=2+k$. Suppose that $k$ is a large nonnegative integer. 
The initial block of $f_q$, that is $x^2y^2$, possesses a square root $r(x)=xy$. However, for any odd $q$ it is impossible to obtain the component $x^q$ as a component of a squared $q$-fps. Thus, there is no $g_q$ for which $f=(g_q)^2$. Observe also that if we continue the computation to determine coefficients $(g_q)_c,\, c\in C_{[2k-1]}$ (according to equations \pref{eq:2} and \pref{eq:4}), we will obtain coefficients of $h$, and the first contradiction in the calculations will appear when we try to compute values $(g_q)_c,\, c\in C_{2k}$. 
\end{exam}


\begin{thebibliography}{10}


\bibitem{WB65} W.G. Brown, \textit{On the existence of square roots in certain rings of power series}, {Math. Ann.} 158 (1965), 82--89.

\bibitem{WB67} W.G. Brown, \textit{On $k^th$ roots in power series rings}, {Math. Ann.} {170} (1967), 327--333.

\bibitem{BTW} W.G. Brown, W. Tutte, \textit{On the enumeration of rooted non-separable planar maps}, Canad. J. Math. 16 (1964), 572--577.



\bibitem{XD} {D. Bugajewski}, {X.-X. Gan}, \textit{Formal multiplicative root series and algorithms of their evaluation}, Comm. Algebra {49} (2021), {3232--3240}.

\bibitem{BBXM} {D. Bugajewski}, {D. Bugajewski}, {X.-X. Gan,} P. Ma\'ckowiak, \textit{On the recursive and explicit form of the general J.C.P. Miller formula with applications}, Adv. Appl. Math. 156 (2024), 102688.

	
\bibitem{CM} M. Cohen, \textit{Roots of formal power series and new theorems on Riordan group elements}, Congr. Numer 233 (2019),  195--204.

\bibitem{CLO} D.A. Cox, J. Little, D. O'Shea, \textit{Ideals, Varietes, and Algorithms}, Springer Undergrad. Texts Math. Technol., Springer, 2015.


\bibitem{Hunk2019} P. Haukkanen, \textit{Formal power series in several variables}, NNTDM 25 (2019), 44--57.

\bibitem{N69} I. Niven, \textit{{Formal power series}}, Am. Math. Mon. {76} (1969), {871--889}.

\bibitem{Samb} B. Sambale, {\it An Invitation to Formal Power Series}, Jahresber. Dtsch. Math. Ver. {125} (2023), 3--69.

\bibitem{S}{R.P. Stanley, \textit{Enumerative Combinatorics}, {Cambridge Univ. Press}, {1997}.}

\bibitem{TW} W. Tutte, \textit{A census of planar triangulations}, {Canad. J. Math.} {14} (1962), {21--38}.

\end{thebibliography}
\end{document}